\documentclass[12pt]{article}
\usepackage[cp866]{inputenc}
\usepackage[english]{babel}
\usepackage{graphicx}
\usepackage{subfig}
\usepackage{extsizes}
\usepackage{epsfig}
\usepackage[numbers,sort&compress]{natbib}
\usepackage{citehack}
\usepackage{float}
\usepackage{amsmath}
\usepackage{amssymb}
\usepackage{multicol}
\usepackage{amscd}
\usepackage{verbatim}
\usepackage{booktabs}
\usepackage{enumitem}
\newcommand{\eps}{\varepsilon}

\newcommand{\vfi}{\varphi}

\newcommand{\R}{\mathbb {R}}

\newcommand{\p}{\partial}

\makeatletter \@addtoreset{equation}{section} \makeatother
\allowdisplaybreaks[1]

\begin{document}
      \title{On change  of slow variables at crossing the separatrices
 }
      \date{}
      \author{Anatoly Neishtadt} 
      \maketitle
    \begin{abstract}
    We consider  general (not necessarily Hamiltonian) perturbations of Hamiltonian systems with one degree of freedom near separatrices of the unperturbed system. We present asymptotic formulas for change of slow variables at evolution across separatrices.
    \end{abstract}
    
\section{Outline of the problem}\label{ss2.1}
We  consider   systems described by  differential equations of the form
 \begin{eqnarray} \label{perturbed}
 \dot q&=&\frac{\partial H}{\partial p}+\eps f_q, \, \dot p=-\frac{\partial H}{\partial q}+\eps f_p, \, \dot z=\eps f_z\, ,
 \\
        H&=&H(p,q,z),\, f_{\alpha}=f_{\alpha}(p,q,z,\varepsilon),\alpha=p,q,z,\, (p,q)\in \R^2,z\in \R^{l-2}\,. \ \nonumber
       \end{eqnarray}
       Here $\varepsilon $  is a small parameter,  $|\varepsilon|\ll 1$.
        For $\varepsilon=0,\,  z={\rm const}$ we have {\it an unperturbed system} for $p,q$, which is a Hamiltonian system with one degree of freedom. The function $H$ is an unperturbed  Hamiltonian.
        For $\eps > 0$ we have {\it a perturbed system},    and  functions $\varepsilon f_{\alpha}$ are {\it  perturbations}.  
        
\begin{figure}
 \begin{center}
            \includegraphics[scale=0.65, angle=0.0]{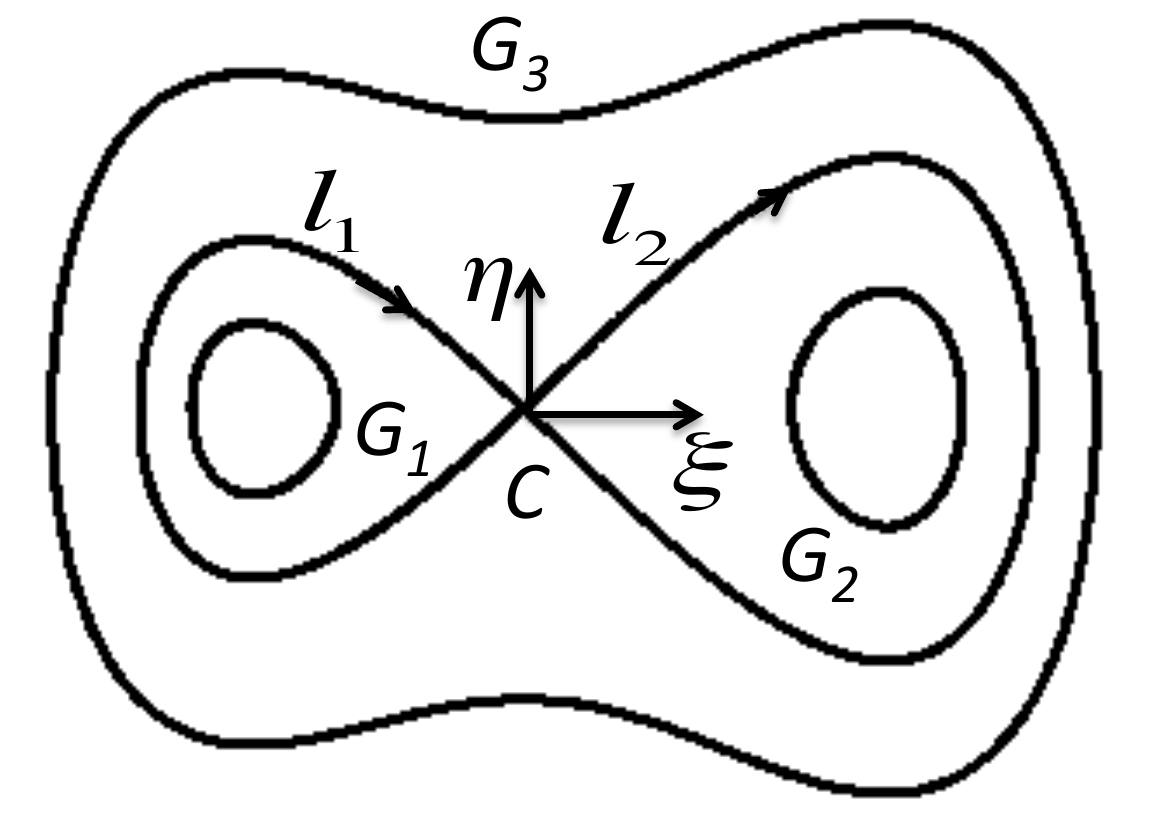}
            \end{center}
           \caption{Phase portrait of the unperturbed system.}
            \label{unperturbed_plane} 
\end{figure}

        It is supposed that there are a saddle point and passing through it separatrices in the phase portrait of the unperturbed system, Fig. \ref{unperturbed_plane}. 
         Under the action of  perturbations 
  the projection of the phase point onto the plane $(p,q)$ crosses a separatrix. 
      
       Separatrices divide the phase plane of the unperturbed systems into domains $G_1(z), G_2(z), G_3(z)$, Fig. \ref{unperturbed_plane}. In each of these domains   it is possible to use  variables $h, \varphi$ instead of  $p,q$, where $h$ is the difference between $H$ and its value at the saddle point, and   $\varphi$ is  ``the angle'' (from the pair ``action-angle''  variables \cite{arn_1} of the unperturbed system). Then for $h,z,\varphi$  we get the perturbed system having the standard form of  system with one rotating phase \cite{bm}: in this system $h,z$ are called {\it slow variables}, $\vfi$ is {\it the rotating phase}. It is a classical result  that the averaged with respect to $\varphi$ system describes the evolution of $h,z$ far from separatrices with accuracy $O(\varepsilon)$ during the time interval of order $1/\varepsilon$ \cite{bm}. For approximate description of evolution of $h,z$ for trajectories that cross the unperturbed separatrices one can  use  the averaged system up to the separatrix and then averaged system with initial conditions on the separatrix  in one of domains in which the trajectory is captured (a certain probability can be assigned to each such continuation). For  majority of  initial conditions this procedure describes the behaviour of  slow variables with accuracy $O(\eps\ln\eps)$ during  time of order $1/\eps$; the measure of the ``bad'' set  of initial conditions, for which this description is not valid, tends to $0$ faster than any given power of $\eps$ as $\eps\to 0$ \cite{Nei_nonlinearity}. One can make one more step of the averaging method and use the same procedure for the second order averaged system (it is shown in \cite{Nei_Okunev_nonlinearity} that solutions of this system indeed arrive to separatrices).     This improves accuracy up to $O(\varepsilon^2)$ for motions far from  separatrices.   However, for motion with separatrix crossing there is no improvement. The reason is that there is a change of order at least $\eps$  of slow variables at crossing  a narrow neighbourhood of separatrices. Because the width of this neighbourhood tends to 0 as $\eps \to 0$, it is reasonable to call this change  {\it  a jump} of slow variables at the separatrix.   In this note we give asymptotic formulas for this jump. Such formulas were first obtained in \cite{Timofeev} for the pendulum in a slowly varying gravitational
field, then in \cite{Cary86, Neishtadt86} for the general case of a Hamiltonian system with one degree of freedom
and slowly varying parameters, in \cite{nei_pmm} for the general case of a slow-fast Hamiltonian system
with two degrees of freedom,  and in \cite{Haberman_B_1990, Haberman_B_1994} for motion in a slowly time-dependent potential with a dissipation.   Jump of slow variables is interpreted as  {\it a jump of an adiabatic invariant} for Hamiltonian systems  \cite{Timofeev, Cary86, Neishtadt86, nei_pmm} and as {\it a time shift} for 
  systems with a   dissipation  \cite{Haberman_B_1990, Haberman_B_1994}.  We consider the case of general perturbed system (\ref{perturbed}).    For derivation of intermediate estimates used in this note see, e.g., \cite{Nei_nonlinearity, Nei_Okunev_nonlinearity}.  
  
       \section{Asymptotic expansions for unperturbed motions near separatrices}
       
       In the phase portrait of the unperturbed system there is a saddle point $C=C(z)$ and passing through it separatrices $l_1=l_1(z), l_2=l_2(z)$. We denote $l_3=l_3(z)=l_1(z)\cup l_2(z)$.  Denote $q_C=q_C(z), p_C=p_C(z)$ coordinates of the point $C$. 
   Denote 
   \begin{equation*}
   \begin{aligned} 
    h_C(z)&=H(p_C(z), q_C(z), z), \ E(p,q,z)=H(p,q,z)-h_C(z).
         \end{aligned}
         \end{equation*}    
                   We assume that $E>0$ in $G_3$, $E<0$ in $G_{1,2}$.

   Denote 
   \begin{equation*}
   \begin{aligned} 
        f_{z,C}(z&)=f_z( p_C(z), q_C(z), z,0),\ F_z(p,q,z)=f(p,q,z,0)-f_{z,C}(z), \\
          f_h(p,q,z)&= \frac{\p E}{\p p}f_p(p,q,z,0)+\frac{\p E}{\p q}f_q(p,q,z,0)+ \frac{\p E}{\p z}f_z(p,q,z,0)         .
         \end{aligned}
         \end{equation*}    
                                      For the period $T$ of the trajectory $E=h$ in  domain $G_i$ we have
          \begin{equation*}
                  \begin{aligned}
                  T=-a_i\ln |h| +b_i +O( h\ln|h|),    \  a_1=a_2=a, a_3=2a, \ b_3=b_1+b_2.
                    \end{aligned}
                  \end{equation*} 
                  
                  Denote
                  \begin{equation*}
                   \oint_{l_i} f_h(p,q,z)dt=-\Theta_i(z), \  \oint_{l_i} F_z(p,q,z)dt=A_i(z).
                  \end{equation*} 
                  Then for integrals along the unperturbed phase trajectory $E=h$ in the domain 
                  $G_i$ we have
                  \begin{equation*}
                  \begin{aligned}
                  \oint_{E=h} f_h(p,q,z)dt&=-\Theta_i(z)+O( h\ln|h|),\\
                  \oint_{E=h} F_z(p,q,z)dt&=A_i(z)+O( h\ln|h|).
                  \end{aligned}
                  \end{equation*} 
                  We assume that $\Theta_1(z)>0, \Theta_2(z)>0$ for all considered values of $z$.
                  
                  \medskip
                  Introduce the coordinate system $C\xi\eta$ as shown in Fig. \ref{unperturbed_plane}. 
                  For initial points on the positive side of the axis $C\eta$ and integrals on the unperturbed phase trajectory $E=h$ (i.e. in $G_3$)    we have
                  \begin{equation*}
                  \begin{aligned}
                  \frac{1}{T}\int_0^T(t-\frac{T}{2})f_h dt&
                  =-\frac{a\ln h (\Theta_2-\Theta_1)/2+ (\Theta_1 b_2-\Theta_2 b_1)/2 +d_3}
                  {-2a\ln h+b_3}    +O(\sqrt h \, ),          \\
                   \frac{1}{T}\int_0^T(t-\frac{T}{2})F_z dt&
                  =-\frac{a\ln h (A_1-A_2)/2- (A_1 b_2-A_2 b_1)/2 +g_3}
                  {-2a\ln h+b_3}    +O(\sqrt h \, ) . 
                                  \end{aligned}
                  \end{equation*} 
                  For initial points on the axis $C\xi$ and integrals on the unperturbed phase trajectory $E=h$    in the domain $G_i, i=1,2$ we have
                  \begin{equation*}
                  \begin{aligned}
                  \frac{1}{T}\int_0^T(t-\frac{T}{2})f_h dt&=-\frac{d_i}
                  {-a\ln |h|+b_1}    +O(\sqrt{ |h|} \, )  ,
\\
                \frac{1}{T}\int_0^T(t-\frac{T}{2})F_zdt&=-\frac{g_i}
                  {-a\ln |h|+b_i}    +O(\sqrt {|h|} \, )  .
                  \end{aligned}
                  \end{equation*}
                  We have $d_3=d_1+d_2, g_3=g_1+g_2$. 
                  
                  In line with the general approach of the averaging method, one can make a  change of variables
\begin{align}
\begin{split} \label{e:varchange}
  h &= \overline h + \varepsilon u_{h, 1}(\overline h, \overline z, \overline \varphi ) + \varepsilon^2 u_{h, 2}(\overline h, \overline z, \overline \varphi ), \\
  z &= \overline{z} + \varepsilon u_{z, 1}(\overline h, \overline z, \overline \varphi) + \varepsilon^2 u_{z, 2}(\overline h, \overline z, \overline \varphi), \\
  \varphi &= \overline \varphi + \varepsilon u_{\varphi, 1}(\overline h, \overline z, \overline \varphi)
\end{split}
\end{align}
that transforms original equations  of motion to the following form:
\begin{align}
\begin{split} \label{e:init-avg}
  \dot{\overline h} &= \varepsilon \overline f_{h, 1}(\overline h, \overline z) + \varepsilon^2 \overline f_{h, 2}(\overline h, \overline z) + \varepsilon^3 \overline f_{h, 3}(\overline h, \overline z, \overline \varphi, \varepsilon), \\
  \dot{\overline z} &= \varepsilon \overline f_{z, 1}(\overline h, \overline z) + \varepsilon^2 \overline f_{z, 2}(\overline h, \overline z) + \varepsilon^3 \overline f_{z, 3}(\overline h, \overline z, \overline \varphi, \varepsilon), \\
  \dot{\overline \varphi} &= \omega(\overline h, \overline z) + \varepsilon \overline f_{\varphi, 1}(\overline h, \overline z) + \varepsilon^2 \overline f_{\varphi, 2}(\overline h, \overline z, \overline \varphi, \varepsilon).
\end{split}
\end{align}
The first order averaged system is obtained by keeping only the first  term in each of these equations.
The second order averaged system  is obtained by neglecting highest order terms in each of these equations. 

One can show that (see \cite{Nei_Okunev_nonlinearity})
\begin{equation*}
u_{h, 1}=  \frac{1}{T}\int_0^T(t-\frac{T}{2})f_h dt, \  u_{z, 1}=  \frac{1}{T}\int_0^T(t-\frac{T}{2})F_z dt.
\end{equation*}
It is convenient to consider evolution using both usual time $t$ and slow time $\eps t$.

\section{Jump of slow variables}
\label{s_jos}

\subsection{General description of motion}

Let a phase point start to move at $t=t_-=0$ (thus $\tau=\tau_-=0$) in the domain $G_3$ at the distance of order 1 from the separatrix.  Denote $h_-, z_-, \varphi_-$ initial values of variables $h, z, \varphi$. Denote $h(t), z(t), \vfi(t)$ solution of the system (\ref{perturbed}) with this initial condition (written in variables $h,z,\vfi$).  The phase point makes rounds close to unperturbed trajectories in $G_3$ while moving closer to the separatrix with each round, approaches the separatrix, crosses the separatrix and continues the motion in domain $G_i$, $i=1$ or $i=2$. Assume, for the sake of being definite, that this is motion in $G_2$. At $t=t_+=K/\eps$  (thus $\tau=\tau_+=K$) the phase point is in $G_2$ at the distance of order 1 form the separatrix. Here $K={\rm const}$. Denote $h_+=h(t_+),
z_+=z(t_+), \vfi_+=\vfi(t_+)$.

\medskip

Denote $\overline h(\tau), \overline z(\tau)$ the solution of the first order averaged system with initial conditions $h_-, z_-$ glued of solutions of  averaged systems for domains $G_3$ and $G_2$ (cf. \cite {Nei_nonlinearity}). Denote $\tau_*$ the moment of the slow time such that $\overline h(\tau_*)=0$ (i.e. $\tau_*$ is the moment of the slow time for the arrival of this solution to the separatrix). Denote $z_*=\overline z(\tau_*)$.

\medskip

Denote $\hat h_-(\tau), \hat z_-(\tau)$ the solution of the second order averaged system with initial, at $\tau=0$, conditions corresponding to $h_-, z_-, \vfi_-$ (i.e., these initial conditions are obtained from $h_-, z_-, \varphi_-$  by transformation (\ref{e:varchange})).  Denote $\hat h_+(\tau), \hat z_+(\tau)$ the solution of the second order averaged system with initial,  at $\tau=\tau_+\,$, conditions corresponding to $h_+, z_+, \vfi_+$. We consider this solution for $\tau\le \tau_+$. Denote $\hat\tau_{*,\mp}$ moments of arrival of these two solutions to the separatrix, $ \hat h_{\mp}(\tau_{*,\mp})=0$.
Denote $\hat z_{*,\mp}=  \hat z_{\mp} (\hat\tau_{*,\mp})$. Denote
\begin{equation}
\Delta\hat \tau_*= \hat\tau_{*,+}-\hat\tau_{*,-}, \ \Delta\hat z_{*}=\hat z_{*,+}-\hat z_{*,-}.
\end{equation}
We will call these values {\it jumps of slow variables at the separatrix}. To estimate these jumps,
we will  consider description of dynamics by the second order averaged system at approaching the separatrix (in $G_3$) and at moving away from the separatrix (in $G_2$).  

For crossing from domain $G_3$ to domain $G_i$, $i=1,2$,  we use also notations  $ \hat \tau_{*,3}=\hat \tau_{*,-},  \hat \tau_{*,i}=\hat \tau_{*,+}, \hat z_{*,3}=\hat z_{*,-},  \hat z_{*,i}=\hat z_{*,+}$.

\medskip
Values $f_{z,C}, \Theta_i, A_i,  a_i, b_i, d_i,g_i$ are taken at $z=z_*$ in all expansions below. 


\subsection{Approaching the separatrix} 
\label{approaching}
Consider motion of the phase point in $G_3$. Projection of the phase point onto $p,q$ plane makes rounds close to unperturbed trajectories while moving closer to   the separatrix with each round.  This projection crosses the ray $C\eta$ on each such round when it moves  close enough to the separatrix. We enumerate $N+1$ moments of time for these intersections starting with the last one: $t_0>t_1>\ldots >t_N>0$. The moment of time $t_N$ is chosen in such a way that for  $0\le t\le t_N$ dynamics of $h, z$ is described with a required (high enough) accuracy  by the second order averaged system, while for $ t_N \le t\le t_0$ expansions near the separatrix can be used for description of motion because the phase point is close enough to the separatrix. 

Denote $\tilde h(t), \tilde z(t), \tilde \vfi(t)$ the result of transformation of solution $h(t), z(t), \vfi(t)$ via formulas (\ref {e:varchange}). Denote
\begin{equation}
\begin{aligned}
h_n=h(t_n), z_n=z(t_n), \hat h_n=\hat h(t_n), \hat z_n=\hat z(t_n), \tilde h_N=\tilde h(t_N), \tilde z_N=\tilde z(t_N).
\end{aligned}
\end{equation}
Denote $ U_{h}(t)= u_{h, 1} + \varepsilon u_{h, 2},\ U_{z}(t)= u_{z, 1} + \varepsilon u_{z, 2} $ where functions $u_{h,i}, u_{z,i}$ are those in (\ref {e:varchange}), and they are calculated at the point $\tilde h(t), \tilde z(t), \tilde \vfi(t)$.
Denote $U_{h,n}=U_{h}(t_n), U_{z,n}=U_{z}(t_n)$.

\medskip
We will use the symbol $\simeq$ in  approximate equalities  without indication of accuracy of the approximation. We have  
$$
z_N= \tilde z_N + \eps U_{z,N}\simeq\hat z_n + \eps U_{z,N}.
$$
Then we have an identity
\begin{equation}
\begin{aligned}
\label{identity_outer}
z_0
&\simeq \hat z_{3,*} +(\hat z|_{h=h_0}- \hat z_{3,*})+\left((\hat z|_{h=h_N}- \hat z|_{h=h_0}) -(z_N-z_0)\right)+(\hat z_N- \hat z|_{h=h_N})    
  +\eps U_{z,N}.
\end{aligned}
\end{equation}
Estimate  terms in this expression separately.
\medskip

a) For $(\hat z|_{h=h_0}- \hat z_{3,*})$.

\medskip
This value is the change of $\hat z$ from the moment of time when $\hat h=0$ 
till the moment of time when $\hat h=h_0$.
In the principal approximation
\begin{equation*}
\begin{aligned}
\dot{\hat z}=\eps(f_{z,C}+\frac{1}{T}A_3), \ \dot{\hat h}=-\eps\frac{1}{T}\Theta_3.
\end{aligned}
\end{equation*}
Hence
\begin{equation*}
\frac{d\hat z}{d \hat h}=-\frac{1}{\Theta_3}(T f_{z,C}+A_3)
\end{equation*}
and
\begin{equation*}
\begin{aligned}
\hat z|_{h=h_0}&- \hat z_{3,*} = -\frac{1}{\Theta_3}\int_0^{h_0}(T f_{z,C}+A_3)dh
= -\frac{f_{z,C}}{\Theta_3}\int_0^{h_0}Tdh-\frac{A_3}{\Theta_3}h_0\\
&=-\frac{f_{z,C}}{\Theta_3}\int_0^{h_0}(-2a\ln h +b_3)dh-\frac{A_3}{\Theta_3}h_0
=-\frac{f_{z,C}}{\Theta_3}\left[-2a(h_0\ln h_0-h_0) +b_3h_0\right]-\frac{A_3}{\Theta_3}h_0.
\end{aligned}
\end{equation*}

b) For $\left((\hat z|_{h=h_N}- \hat z|_{h=h_0})-(z_N-z_0)\right)$.

\medskip
To calculate this term one can consider motion round by round, calculate differences between changes of $\hat z$ and $z$  on each round, and sum up these differences.   For  changes of $h, z$ one can use
\begin{equation*}
\begin{aligned}
&h_{n+1}-h_n \simeq  \eps\Theta_3,\\
&z_{n+1}-z_n \simeq  -\eps f_{z, C}\left(-\frac{a}{2}\ln h_n- a\ln(h_n+\eps \Theta_1) -\frac{a}{2}\ln h_{n+1} +b_3 \right) -\eps A_3.
\end{aligned}
\end{equation*} 
The change of $\hat z$ is calculated as
\begin{equation*}
\hat z|_{h=h_{n+1}}- \hat z|_{h=h_n}\simeq-\frac{f_{z, C}}{\Theta_3}\int_{h_n}^{h_{n+1}}(-2a\ln h +b_3)dh-\eps A_3.
\end{equation*}
Thus
\begin{equation*}
\begin{aligned}
&(z_{n+1}-z_n)-(\hat z|_{h=h_{n+1}}- \hat z|_{h=h_n})\\
&\simeq  a
 \frac {f_{z, C}}{\Theta_3}
 \left[\int_{h_n}^{h_{n+1}}(-2\ln h )dh-\eps\Theta_3\left(-\frac{1}{2}\ln h_n- a\ln(h_n+\eps\Theta_1) 
-\frac{1}{2}\ln h_{n+1}\right) \right]
\end{aligned}
\end{equation*}
The expression in the square brackets is related to calculation of the integral of $-\ln h$ by the trapezoidal method like in \cite {Neishtadt86}. Thus, we can directly use the expression for change of an adiabatic invariant from  \cite {Neishtadt86}. This gives
\begin{equation*}
\begin{aligned}
&(\hat z|_{h=h_N}- \hat z|_{h=h_0})-(z_N-z_0)
\\
&\simeq 2\eps a f_{z,C}
\left[-\frac{1}{2} \ln\frac{2\pi}{\Gamma(\xi_3)\Gamma(\xi_3+\theta_{13})} +\xi_3+\left(-\xi_3+\frac{1}{2}\theta_{23}\right)\ln \xi_3\right]\\
&+\eps \frac{1}{2}af_{z,C}(\theta_{23}-\theta_{13})(\ln h_N
- \ln h_0).
\end{aligned}
\end{equation*}
Here $\Gamma(\cdot)$ is the gamma function, $\xi_3=h_0/\Theta_3$,   $\theta_{ij}=\Theta_i/\Theta_j$.
\newpage
c) For $(\hat z_N- \hat z|_{h=h_N}) $.

\medskip
We have $h(t_N)=h_N$. Denote $\hat t_N$ the moment of time such that $\hat h(\hat t_N)=h_N$. Find $\hat t_N -t_N$.

We have
\begin{equation*}
\hat h(\hat t_N)=h_N=h(t_N)=\tilde h(t_N)+\eps U_{h,N}\simeq \hat h(t_N)+\eps U_{h,N}.
\end{equation*}
Thus
\begin{equation*}
\hat t_N -t_N\simeq-\frac{1}{\Theta_3}TU_{h,N}.
\end{equation*}
Value of $T$ is calculated at $h=h_N$. 
Then
\begin{equation*}
\begin{aligned}
\hat z_N&- \hat z|_{h=h_N}=\hat z(t_N)-\hat z(\hat t_N)\simeq \eps\left (f_{z,C}+\frac{1}{T}A_3\right)(t_N-\hat t_N)
\simeq \eps\left (f_{z,C}+\frac{1}{T}A_3\right)\frac{1}{\Theta_3}TU_{h,N}\\
&\simeq-\eps\frac {f_{z,C}}{\Theta_3}
\left(\frac{a}{2}(\Theta_2-\Theta_1)\ln h_N + (\Theta_1 b_2-\Theta_2 b_1)/2 +d_3\right)
 +\eps \frac{1}{4}A_3(\theta_{23}-\theta_{13}).
\end{aligned}
\end{equation*}

d) For $\eps U_{z,N}$.

\medskip
We have  $\eps U_{z,N}\simeq \eps(A_1-A_2)/4$.

\medskip
Combining results of a) - d) we get from identity (\ref{identity_outer})
\begin{equation}
\begin{aligned}
\label{combined_outer}
z_0
&\simeq \hat z_{3,*} -\frac{f_{z,C}}{\Theta_3}\left[-2a(h_0\ln h_0-h_0) +b_3h_0\right]-\frac{A_3}{\Theta_3}h_0\\
&+2\eps a f_{z,C}
\left[-\frac{1}{2} \ln\frac{2\pi}{\Gamma(\xi_3)\Gamma(\xi_3+\theta_{13})} +\xi_3+\left(-\xi_3+\frac{1}{2}\theta_{23}\right)\ln \xi_3\right]\\
&-\eps \frac{1}{2}af_{z,C}(\theta_{23}-\theta_{13})
 \ln h_0\\
&-\eps\frac{1}{2}f_{z,C}
\left( (\theta_{13} b_2-\theta_{23} b_1) +2\frac{d_3}{\Theta_3}\right)
 +\eps \frac{1}{4}A_3(\theta_{23}-\theta_{13}) \\  
 &+ \eps(A_1-A_2)/4\\
 &=\hat z_{3,*} -\frac{f_{z,C}}{\Theta_3}\left[-2a h_0\ln (\eps\Theta_3) +b_3h_0\right]-\frac{A_3}{\Theta_3}h_0\\
&+2\eps a f_{z,C}
\left[-\frac{1}{2} \ln\frac{2\pi}{\Gamma(\xi_3)\Gamma(\xi_3+\theta_{13})} +\frac{1}{2}\theta_{23}\ln \xi_3\right]\\
&-\eps \frac{1}{2}af_{z,C}(\theta_{23}-\theta_{13})
 \ln h_0\\
&-\eps\frac{1}{2}f_{z,C}
\left( (\theta_{13} b_2-\theta_{23} b_1) +2\frac{d_3}{\Theta_3}\right)
 +\eps \frac{1}{4}A_3(\theta_{23}-\theta_{13}) \\  
 &+ \eps(A_1-A_2)/4.
\end{aligned}
\end{equation}

\subsection{Passage through the separatrix}
\label{passage}
We assume that  $ k\sqrt \eps  \le   \xi_3 \le \theta_{23}-  k\sqrt \eps$, where $k$ is a large enough constant. Then for $t>t_0$ the phase point  makes a round close to the separatrix $l_2$ and arrives to the ray $C\xi$ in $G_2$  (Fig. \ref{unperturbed_plane}) at some moment of time $t'_0 =t_0+O(\ln \eps)$. Denote $h_0'=h(t'_0), z_0'=z(t'_0)$. We have
\begin{equation}
\label{combined_passage}
h_0'\simeq h_0-\eps \Theta_2,\  z_0'\simeq z_0+\eps f_{z,C}\left[  -\frac{a}{2}\ln h_0  -\frac{a}{2}\ln (-h_0')  +b_2 \right]+\eps A_2.
\end{equation}
Denote $\xi_2 =(-h_0')/(\eps\Theta_2)\simeq (\eps \Theta_2-  h_0)/(\eps\Theta_2)=(\eps \Theta_2-  \eps \Theta_3 \xi_3)/(\eps\Theta_2)=1-(\Theta_3/\Theta_2)\xi_3$.
Thus $\xi_3\simeq\theta_{23}(1-\xi_2)$.

\subsection{Moving away from the separatrix}
\label{departure}
For $t >t_0'$ the projection of the phase point onto $p,q$ plane makes rounds close to unperturbed trajectories while moving farther away from  the separatrix with each round.  This projection crosses the ray $C\xi$ in $G_2$ on each such round while it moves  close enough to the separatrix. We enumerate $N+1$ moments of time for these intersections starting with the first one: $t_0'<t_1'<\ldots <t_N'<K/\eps$. The moment of time $t_N'$ is chosen in such a way that for 
$t_N'\le t\le K/\eps$
changes of $h, z$ are described with a required (high enough) accuracy  by the second order averaged system, while for  $t_0' \le t\le t_N'$   expansions near the separatrix can be used for description of motion because the phase point is close enough to the separatrix. Calculations here are similar   to those for approaching the separatrix in Section \ref{approaching}. In what follows, we omit `prime' in notation for moments of time and use for variables  $h,z$  the same notation as in Section \ref{approaching}, except of $h_0', z_0'$.  

\medskip
We have
$$
z_N= \tilde z_N + \eps U_{z,N}\simeq\hat z_n + \eps U_{z,N}.
$$
Then we have an identity
\begin{equation}
\begin{aligned}
\label{identity_inner}
z_0'
&\simeq \hat z_{2,*} +(\hat z|_{h=h_0'}- \hat z_{2,*})+\left((\hat z|_{h=h_N}- \hat z|_{h=h_0'}) -(z_N-z_0')\right)+(\hat z_N- \hat z|_{h=h_N})    
  +\eps U_{z,N}.
\end{aligned}
\end{equation}
Estimate  terms in this expression separately.
\medskip

a) For $(\hat z|_{h=h_0}- \hat z_{2,*})$.
\medskip

Similarly to  Section \ref{approaching} we get 
\begin{equation}
\hat z|_{h=h_0'}- \hat z_{2,*}\simeq -\frac{f_{z,C}}{\Theta_2}\left[-a(h_0'\ln |h_0'|-h_0') +b_2h_0'\right]-\frac{A_2}{\Theta_2}h_0'.
\end{equation}

b) For $\left((\hat z|_{h=h_N}- \hat z|_{h=h_0'}) -(z_N-z_0')\right)$.

\medskip
Similarly to Section \ref{approaching} we can use result of \cite {Neishtadt86}. This gives 
\begin{equation}
(\hat z|_{h=h_N}- \hat z|_{h=h_0'}) -(z_N-z_0')\simeq -\eps a f_{z,C}\left[- \ln \frac{\sqrt{2\pi}}{\Gamma(\xi_2)} 
+\xi_2+(\frac{1}{2}-\xi_2)\ln \xi_2\right].
\end{equation}

c) For $(\hat z_N- \hat z|_{h=h_N})$

\medskip
Similarly to Section \ref{approaching} we get
\begin{equation}
\hat z_N- \hat z|_{h=h_N}\simeq-\eps\frac{f_{z,C}}{\Theta_2}d_2.
\end{equation}

d) For $\eps U_{z,N}$.
\medskip

We get $\eps U_{z,N}\simeq 0$.

\medskip
Combining results of a) - d) we get from identity (\ref{identity_inner})

\begin{equation}
\begin{aligned}
\label{combined_inner}
z_0'
&\simeq \hat z_{2,*} -\frac{f_{z,C}}{\Theta_2}\left[-a(h_0'\ln |h_0'|-h_0') +b_2h_0'\right]-\frac{A_2}{\Theta_2}h_0'\\
&-\eps a f_{z,C}\left[- \ln \frac{\sqrt{2\pi}}{\Gamma(\xi_2)} 
+\xi_2+(\frac{1}{2}-\xi_2)\ln \xi_2\right]-\eps\frac{f_{z,C}}{\Theta_2}d_2\\
&=\hat z_{2,*} -\eps{f_{z,C}} \left[a(\xi_2\ln (\eps\Theta_2\xi_2)-\xi_2) - b_2\xi_2\right]+\eps {A_2}\xi_2\\
&-\eps a f_{z,C}\left[- \ln \frac{\sqrt{2\pi}}{\Gamma(\xi_2)} 
+\xi_2+(\frac{1}{2}-\xi_2)\ln \xi_2\right]-\eps\frac{f_{z,C}}{\Theta_2}d_2\\
&=\hat z_{2,*} -\eps{f_{z,C}} \left[a\xi_2\ln (\eps\Theta_2) - b_2\xi_2\right]+\eps {A_2}\xi_2\\
&-\eps a f_{z,C}\left[- \ln \frac{\sqrt{2\pi}}{\Gamma(\xi_2)} 
+\frac{1}{2}\ln \xi_2\right]-\eps\frac{f_{z,C}}{\Theta_2}d_2.
\end{aligned}
\end{equation}

\subsection{Formula for jump of slow variables}
\label{jump_1}

Combining results of Sections \ref{approaching}, \ref{passage}, \ref{departure}  (formulas (\ref {combined_outer}),
(\ref {combined_passage}) and (\ref {combined_inner})  ) we get
\begin{equation}
\begin{aligned}
&\hat z_{2,*} -\eps{f_{z,C}} \left[a\xi_2\ln (\eps\Theta_2) - b_2\xi_2\right]+\eps {A_2}\xi_2\\
&-\eps a f_{z,C}\left[- \ln \frac{\sqrt{2\pi}}{\Gamma(\xi_2)} 
+\frac{1}{2}\ln \xi_2\right]-\eps\frac{f_{z,C}}{\Theta_2}d_2\\
&\simeq \hat z_{3,*} -\frac{f_{z,C}}{\Theta_3}\left[-2a h_0\ln (\eps\Theta_3) +b_3h_0\right]-\frac{A_3}{\Theta_3}h_0\\
&+2\eps a f_{z,C}
\left[-\frac{1}{2} \ln\frac{2\pi}{\Gamma(\xi_3)\Gamma(\xi_3+\theta_{13})} +\frac{1}{2}\theta_{23}\ln \xi_3\right]\\
&-\eps \frac{1}{2}af_{z,C}(\theta_{23}-\theta_{13})
 \ln h_0\\
&-\eps\frac{1}{2}f_{z,C}
\left( (\theta_{13} b_2-\theta_{23} b_1) +2\frac{d_3}{\Theta_3}\right)
 +\eps \frac{1}{4}A_3(\theta_{23}-\theta_{13}) \\  
 &+ \eps(A_1-A_2)/4\\
 &+\eps f_{z,C}\left[  -\frac{a}{2}\ln h_0  -\frac{a}{2}\ln (-h_0')  +b_2 \right]+\eps A_2
\end{aligned}
\end{equation}
Therefore
\begin{equation}
\begin{aligned}
\label{for_jump_0}
&\Delta\hat z_{*}=\hat z_{*,+}-\hat z_{*,-}=\hat z_{2,*}- \hat z_{3,*}\\
&\simeq \eps{f_{z,C}} \left[a\xi_2\ln (\eps\Theta_2) - b_2\xi_2\right]-\eps {A_2}\xi_2\\
&\eps a f_{z,C}\left[- \ln \frac{\sqrt{2\pi}}{\Gamma(\xi_2)} 
+\frac{1}{2}\ln \xi_2\right]+\eps\frac{f_{z,C}}{\Theta_2}d_2\\
&-\frac{f_{z,C}}{\Theta_3}\left[-2a h_0\ln (\eps\Theta_3) +b_3h_0\right]-\frac{A_3}{\Theta_3}h_0\\
&+2\eps a f_{z,C}
\left[-\frac{1}{2} \ln\frac{2\pi}{\Gamma(\xi_3)\Gamma(\xi_3+\theta_{13})} +\frac{1}{2}\theta_{23}\ln \xi_3\right]\\
&-\eps \frac{1}{2}af_{z,C}(\theta_{23}-\theta_{13})
 \ln h_0\\
&-\eps\frac{1}{2}f_{z,C}
\left( (\theta_{13} b_2-\theta_{23} b_1) +2\frac{d_3}{\Theta_3}\right)
 +\eps \frac{1}{4}A_3(\theta_{23}-\theta_{13}) \\  
 &+ \eps(A_1-A_2)/4\\
 &+\eps f_{z,C}\left[  -\frac{a}{2}\ln h_0  -\frac{a}{2}\ln (-h_0')  +b_2 \right]+\eps A_2.
 \end{aligned}
\end{equation}
For passage form $G_3$ to $G_1$ we would have  relation  (\ref{for_jump_0}) with replacement of index `2' by index `1' . The final result in the form for passage from $G_3$ to $G_i$ where $i=1 \ {\rm or} \ 2$  is simplified to
\begin{equation}
\begin{aligned}
\label{z_3i}
&\Delta\hat z_{*}= \hat z_{i,*}- \hat z_{3,*}
\simeq\eps{f_{z,C}}a(\xi_i-\frac{1}{2})(\ln(\eps\Theta_i) -2\theta_{i3}\ln(\eps\Theta_3))\\
& -\eps a f_{z,C}\ln{\frac{(2\pi)^{3/2}}{\Gamma(\xi_i)\Gamma(\theta_{i3}(1-\xi_i)\Gamma(1-\theta_{i3}\xi_i)} }\\
&   -\eps f_{z,C} (\xi_i-\frac{1}{2} )(b_i- \theta_{i3}b_3) -\eps(\xi_i-\frac{1}{2} )(A_i- \theta_{i3}A_3)\\
& +\eps\frac {f_{z,C}}{\Theta_i}\left({d_i} -  \theta_{i3}{d_3}  \right).
   \end{aligned}
\end{equation}

This formula is the main result of the current note. In a similar way one can write  formulas  for jumps of slow variables due to other passages between domains $G_j$ that occur for other signs of values $\Theta_j, j=1,2,3$.

\medskip
Value $\xi_i$ is called  {\it a crossing parameter} or {\it a pseudo-phase}. Asymptotic formulas for the pseudo-phase were obtained in \cite{Cary_Skodje} for Hamiltonian systems with one degree of freedom and slow
time dependence, in \cite{Neishtadt_Vasiliev}  for slow-fast Hamiltonian systems with one degree of freedom
corresponding to fast motion,  in \cite{Haberman_B_1990, Haberman_B_1994} for motion in a slowly time-dependent potential with a dissipation, and in \cite{Nei_Okunev_nonlinearity} for a general perturbed system of form (\ref {perturbed}).

\medskip

{\bf Remark.} We do not indicate  accuracy of  formula (\ref {z_3i}). One can see that  terms $\sim \eps/\ln h_N$ are neglected in some intermediate relations. However, because the final result should not depend on  choice of $h_N$, the accuracy of the final formula should be much better. For Hamiltonian perturbations the accuracy of the final formula is 
$ O(\eps^{3/2}(|\ln \eps|+(1-\xi_i)^{-1}))$  \cite{Neishtadt86,nei_pmm}.

\subsection{Shift of slow time}
\label{jump_t}

The slow time $\tau$ can be considered as a particular slow variable, $\dot \tau =\eps$. The formula for jump (or {\it sift}) of slow time for  passage from $ G_3$ to $G_i$,  $i=1\ {\rm or }\ 2$, is a particular case   of (\ref{z_3i}) with $f_{z,C}=1$, $A_j=0, j=1,2,3$. Thus we get
\begin{equation}
\begin{aligned}
\label{tau_3i}
 &\hat \tau_{i,*}- \hat \tau_{3,*}
\simeq\eps a(\xi_i-\frac{1}{2})(\ln(\eps\Theta_i) -2\theta_{i3}\ln(\eps\Theta_3))\\
& -\eps a \ln{\frac{(2\pi)^{3/2}}{\Gamma(\xi_i)\Gamma(\theta_{i3}(1-\xi_i)\Gamma(1-\theta_{i3}\xi_i)} }\\
&   -\eps (\xi_i-\frac{1}{2} )(b_i- \theta_{i3}b_3) 
 + \frac{\eps}{\Theta_i}\left({d_i}
  -  \theta_{i3}{d_3}  \right).
   \end{aligned}
\end{equation}

\section{Jump of adiabatic invariant}
In this Section we derive formulas for jumps of adiabatic invariants in Hamiltonian systems from  obtained formulas for jumps of slow variables. 

\subsection{Time-dependent Hamiltonian system}
Let  system (\ref{perturbed}) be a  Hamiltonian system with the Hamiltonian $H = H(p,q.\tau)$, $\tau= \eps t$. Denote  $S_i(\tau)$ area of the domain $G_i$, $i=1,2$. Denote
 $S_3(\tau)= S_1(\tau)\cup S_2(\tau)$. Then $\Theta_j=dS_j/d\tau$, $j=1,2,3$.
 Consider motion with passage from $G_3$ to $G_i$, $i=1 \ {\rm  or} \ 2$, as in Section \ref{s_jos}.  Let $J_-$ and $J_+$ be the initial (at $t=0$, in $G_3$) and final (at $t=K/\eps$, in $G_i$) values of the improved adiabatic invariant. (For the definition of the  improved adiabatic invariant and related formulas see, e.g., \cite{Neishtadt86}). 
Then $S_3(\hat \tau_{3,*})\simeq 2\pi J_-$, $ S_i(\hat \tau_{i,*})\simeq 2\pi J_+$. 
We get
\begin{equation}
\begin{aligned}
2\pi J_+\simeq  S_i(\hat \tau_{i,*})=S_i(\hat \tau_{3,*}  + \hat \tau_{i,*}-\hat \tau_{3,*})  \simeq  S_i(\hat \tau_{3,*}) +\Theta_i(\hat \tau_{i,*} - \hat \tau_{3,*}).
\end{aligned}
\end{equation}
Substitute $(\hat \tau_{i,*}- \hat \tau_{3,*})$ from (\ref{tau_3i}). We get
\begin{equation}
\begin{aligned}
&2\pi J_+\simeq  S_i(\hat \tau_{3,*})+\eps a\Theta_i (\xi_i-\frac{1}{2})(\ln(\eps\Theta_i) -2\theta_{i3}\ln(\eps\Theta_3))\\
& -\eps a \Theta_i \ln{\frac{(2\pi)^{3/2}}{\Gamma(\xi_i)\Gamma(\theta_{i3}(1-\xi_i)\Gamma(1-\theta_{i3}\xi_i)} }\\
&   -\Theta_i (\xi_i-\frac{1}{2} )(b_i- \theta_{i3}b_3) 
 + {\eps}\left({d_i}
  -  \theta_{i3}{d_3}  \right)
\end{aligned}
\end{equation}
as in \cite{Cary86, Neishtadt86}. One can replace $S_i(\hat \tau_{3,*})$ with $S_i(\tau_{*})+\theta_{i3} ( 2\pi J_- -S_3(\tau_{*}))$  here.

\subsection{Slow-fast Hamiltonian system}

Let  system (\ref{perturbed}) be a slow-fast Hamiltonian system. The Hamiltonian is $H(p,q,y,x)$ with pairs of conjugate variables $(p, q)$ and $(y, \eps^{-1} x)$.     Equations of motion are 
\begin{eqnarray} \label{slow-fast}
 \dot q=\frac{\partial H}{\partial p} , \ \dot p=-\frac{\partial H}{\partial q}, \
   \dot x=\eps \frac{\partial H}{\partial y }, \  \dot y=-\eps \frac{\partial H}{\partial x }.
       \end{eqnarray}
       Thus, $z=(y,x)$,  $f_{z,C}=( -{\partial h_C(y,x)}/{\partial x }, {\partial h_C(y,x)}/{\partial x })$.
       Denote  $S_i(z)=S_i(y,x)$ area of domain $G_i$, $i=1,2$. Denote $S_3(z)=S_1(z)\cup S_2(z)$. Then $\Theta_j=\{S_j,h_c \}, j=1,2,3$, where $\{\cdot, \cdot\}$ is the Poisson bracket with respect to variables $(y,x)$,
   $ \{a, b\}= a'_x b'_y- a'_y b'_x$  
        (see \cite{nei_pmm}) and 
\begin{equation}
A_j=\left(-\oint_{l_j} \frac{\p (H-h_c) }{\p x } dt, \oint_{l_j} \frac{\p (H-h_c) }{\p y } dt\right)=\left( \frac{\p S_j }{\p x }, -\frac{\p S_j }{\p y }\right). 
\end{equation}  
(cf.     \cite{nei_pmm}).

Consider motion with passage from $G_3$ to $G_i$, $i=1 \ {\rm  or} \ 2$, as in Section \ref{s_jos}.  Let $J_-$ and $J_+$ be the initial (at $t=0$, in $G_3$) and final (at $t=K/\eps$, in $G_i$) values of the improved adiabatic invariant. (For the definition of the  improved adiabatic invariant and related formulas see, e.g., \cite{nei_pmm}). 
Then $S_3(\hat z_{3,*})\simeq 2\pi J_-$, $S_i(\hat z_{i,*})\simeq 2 \pi J_+$. 
Then we get
\begin{equation}
\begin{aligned}
2\pi J_+\simeq  S_i(\hat z_{i,*})=S_i(\hat z_{3,*} +\hat z_{i,*}-\hat z_{3,*})\simeq  S_i(\hat z_{3,*}) +({\rm grad} \, S_i \cdot( \hat z_{i,*}- \hat z_{3,*})   ).
\end{aligned}
\end{equation}
Here $(\phantom {*} \cdot \phantom{*})$ is the standard scalar product.
Substitute  $(\hat z_{i,*}- \hat z_{3,*})$ from (\ref{z_3i}) and note that
\begin{equation}
\begin{aligned}
({\rm grad} \, S_i \cdot f_{z,C})=\Theta_i, \
({\rm grad} \, S_i \cdot A_i)=\{S_i,S_i\}=0, \ ({\rm grad} \, S_i \cdot A_3)=-\{S_i,S_3\}.
\end{aligned}
\end{equation}
We get
\begin{equation}
\begin{aligned}
\label{J_sf}
&2\pi J_+
\simeq  S_i(\hat z_{3,*})+
\eps \Theta_i a(\xi_i-\frac{1}{2})(\ln(\eps\Theta_i) -2\theta_{i3}\ln(\eps\Theta_3))\\
& -\eps a \Theta_i  \ln{\frac{(2\pi)^{3/2}}{\Gamma(\xi_i)\Gamma(\theta_{i3}(1-\xi_i)\Gamma(1-\theta_{i3}\xi_i)} }\\
&   -\eps \Theta_i  (\xi_i-\frac{1}{2} )(b_i- \theta_{i3}b_3) -\eps \theta_{i3}(\xi_i-\frac{1}{2} ) \{S_i,S_3\}\\
& +\eps\left({d_i} -  \theta_{i3}{d_3}  \right).
\end{aligned}
\end{equation}
For systems with two degrees of freedom  one can approximately calculate $S_i(\hat z_{3,*})$ via initial value of the improved adiabatic invariant and solution of the first order averaged system. Consider motion in the energy level $H=h$. 
Relations
   \begin{equation}  
   S_3(\hat z_{3,*})\simeq 2\pi J_-, \quad   h_C(\hat z_{3,*})\simeq h, \quad h_C(z_*)=h
   \end {equation} 
   imply that
 \begin{equation}
\begin{aligned}
\label{eq_h_C}
({\rm grad} \, S_3\cdot (\hat z_{3,*}-z_*))\simeq 2\pi J_- - S_3(z_*), \quad ({\rm grad} \, h_C \cdot (\hat z_{3,*}-z_*))\simeq 0.
\end{aligned}
\end{equation}
We have
 \begin{equation}
 \label{Si_appr}
S_i(\hat z_{3,*})=S_i(z_*+\hat z_{3,*}-z_*)\simeq S_i(z_{*})+({\rm grad} \, S_i\cdot (\hat z_{3,*}-z_*)).
\end{equation}  
Solve equations (\ref {eq_h_C}) for $(\hat z_{3,*}-z_*)$ and substitute the result to (\ref {Si_appr}). We get
\begin{equation}
S_i(\hat z_{3,*})\simeq S_i(z_{*})+\theta_{i3} (2\pi J_-- S_3(z_*)).
\end{equation}
Substitution of this relation to (\ref{J_sf}) gives an expression for jump of the adiabatic invariant in  \cite{nei_pmm}. 

\section{Conclusions}
The main result of this note is the asymptotic  formula (\ref{z_3i}) for change of slow variables at evolution across separatrices in systems of form (\ref{perturbed}). Together with formula  for phase change in such systems \cite{Nei_Okunev_nonlinearity} this gives a rather complete description of dynamics with separatrix crossings in the considered class of systems.

\vskip 5mm

\noindent Anatoly Neishtadt

\noindent {\small Department of Mathematical Sciences}

\noindent {\small Loughborough University, Loughborough LE11 3TU, United Kingdom}

\noindent {\footnotesize{E-mail: a.neishtadt@lboro.ac.uk}}


\newpage

\begin{thebibliography}{99}

  \bibitem{arn_1} {Arnold V.~I. \/} Mathematical Methods of Classical Mechanics: Graduate Texts in Mathematics 60. Springer-Verlag, New York (1978), x+462 pp.






  \bibitem{bm} { Bogolyubov N.~N.,  Mitropolskij Yu.~A.\/} Asymptotic Methods in the Theory of Non-Linear Oscillations. Hindustan Publ. Corp., Delhi; Gordon and Breach Sci.
Publ., New York (1961), x+537 pp.

\bibitem{Haberman_B_1990} Bourland F. J., Haberman R. Separatrix crossing: time-invariant potentials with dissipation. {\em SIAM J. Appl. Math.},  {\bf 50}, 6,  1716--1744 (1990)


\bibitem{Haberman_B_1994} Bourland F. J., Haberman R.  Connection across a separatrix with dissipation. {\it Stud. Appl. Math.},  {\bf 91},  95--124 (1994) 



\bibitem{Cary86}  
{Cary J. R.}, {Escande D. F.}, {Tennyson J. L.} 
\newblock {Adiabatic-invariant change due to separatrix crossing}.
\newblock {\em Physical Review A}, {\bf 34},  4256--4275 (1986)

\bibitem{Cary_Skodje}  Cary J. R.,  Skodje R. T.  Phase change between separatrix crossings. {\it Physica D},  {\bf 36}, 3,  287--316 (1989)

\bibitem{Neishtadt86}
{Neishtadt A. I. }  
\newblock {Change of an adiabatic invariant at a separatrix}.
\newblock {\em Soviet Journal of Plasma Physics}, {\bf 12}, 568--573 (1986)


  \bibitem{nei_pmm}{Neishtadt A. I.} 
  On the change in the adiabatic invariant on crossing a separatrix
in systems with two degrees of freedom. {\it  J. Appl. Math. Mech.},  {\bf 51},  5,  586--592 (1987)

\bibitem{Neishtadt_Vasiliev} Neishtadt A. I.,  Vasiliev A. A. Phase change between separatrix crossings in slow-fast
Hamiltonian systems.{\it  Nonlinearity},  {\bf 18}, 3, 1393--1406 (2005)












\bibitem {Nei_nonlinearity} Neishtadt A. I.   Averaging method for systems with separatrix crossing. {\it Nonlinearity}, {\bf  30}, 7,  2871--2917 (2017)



\bibitem {Nei_Okunev_nonlinearity} Neishtadt A. I., Okunev A.V.   Phase change and order 2 averaging for one-frequency systems with separatrix crossing. {\it Nonlinearity},  {\bf 35},  8,  4469--4516 (2022) 

 \bibitem{Timofeev}   {Timofeev A. V. }  On the constancy of an adiabatic invariant when the nature
of the motion changes.  {\it Sov. Phys., JETP}, {\bf 48},  656--659 (1978)

\end{thebibliography}
\end{document}